\documentclass[12pt,reqno]{amsart}
\usepackage{amssymb}

\numberwithin{equation}{section}
\theoremstyle{plain}
\newtheorem{theorem}{Theorem}[section]

\newtheorem{proposition}[theorem]{Proposition}
\newtheorem{corollary}[theorem]{Corollary}
\newtheorem{lemma}[theorem]{Lemma}

\theoremstyle{definition}
\newtheorem{remark}[theorem]{Remark}

\def \N {\mathbb{N}}
\def \R {\mathbb{R}}

\def \E {\mathbb{E}}
\def \P {\mathbb{P}}

\def \e {\varepsilon}
\def \d {\delta}
\def \D {\Delta}

\def \s {\sigma}

\def \< {\langle}
\def \> {\rangle}

\def \rank {{\rm rank}}

\def \Col {{\rm Col}}

\newcommand{\norm}[1]{\left \| #1 \right \|}
\newcommand{\sumi}{\sum_{i=1}^d}
\newcommand{\tens}[1]{#1 \otimes #1}

\begin{document}
\title [Sampling from large matrices]
  {Sampling from large matrices: an approach through
  geometric functional analysis }

\author{Mark Rudelson}
\address{Department of Mathematics, University of Missouri, Columbia, MO 65211, U.S.A.}
\email{rudelson@math.missouri.edu}

\author{Roman Vershynin}
\address{Department of Mathematics, University of California, Davis, CA 95616, U.S.A.}
\email{vershynin@math.ucdavis.edu}

\thanks{{\it The ACM Computing Classification:} {\bf G.1.3} Numerical Linear Algebra,
{\bf F.2.1} Numerical Algorithms and Problems}

\thanks{The first author is partially supported by the NSF grant DMS 0245380.
  The second author is an Alfred P. Sloan Research Fellow.
  He is also partially supported by the NSF grant DMS 0401032. }

\date{\today}

\begin{abstract}
We study random submatrices of a large matrix $A$.
We show how to approximately compute $A$ from its random submatrix
of the smallest possible size $O(r \log r)$ with a small error in
the spectral norm, where $r = \|A\|_F^2 / \|A\|_2^2$
is the numerical rank of $A$.
The numerical rank is always bounded by, and is a stable relaxation of,
the rank of $A$.
This yields an asymptotically optimal guarantee in an algorithm
for computing {\em low-rank approximations} of $A$.
We also prove asymptotically optimal estimates on the {\em spectral norm}
and the {\em cut-norm} of random submatrices of $A$.
The result for the cut-norm yields a slight improvement on the best known
sample complexity for an approximation algorithm for MAX-2CSP problems.
We use methods of Probability in Banach spaces, in particular
the law of large numbers for operator-valued random variables.
\end{abstract}

\maketitle

\section{Introduction}

This paper studies random submatrices of a large matrix $A$. The
study of random submatrices spans several decades and is related
to diverse areas of mathematics and computer science. Two main
reasons for the interest in random submatrices are:
\begin{enumerate}
\item one can learn properties of $A$ from the properties of its random submatrices;
\item properties of $A$ may improve by passing to
  its random submatrices.
\end{enumerate}
We address both aspects of random submatrices in this paper.
We show how to approximate $A$ by its random submatrix in
the spectral norm, and we compute the asymptotics of
the spectral and the cut norms of random submatrices.
This yields improvements upon known algorithms for
computing low rank approximations, Singular Value Decompositions,
and approximations to MAX-2CSP problems.

\subsection{The spectral norm: low rank approximations and SVD}

Can one approximate $A$ by only knowing a random submatrix of $A$
of a fixed size? If so, what is the sample complexity,
the minimal size of a submatrix which yields a good approximation
with a small error in some natural norm, and with high probability?

This problem belongs to a class of problems the Statistical
Learning Theory is concerned with. These problems inevitably bear
the assumption that the the object to be learned belongs to a
relatively small ``target'' class. To be able to learn $A$ from a
matrix of small size thus of small rank, we have to assume that
$A$ itself has {\em small rank}--or can be approximated by an
(unknown) matrix of a small rank. We thus strive to find a low
rank approximation to a matrix $A$, whenever such an approximation
exists, from only knowing a small random submatrix of $A$.

Solving this problem is essential for development of fast
Monte-Carlo algorithms for computations on large matrices. An
extremely large matrix $A$ -- say, of the order of $10^5
\times 10^5$ -- is impossible to upload into the Random Access
Memory (RAM) of a computer; it is instead stored in an external
memory. On the other hand, sampling a small submatrix of $A$,
storing it in RAM and computing its small rank approximation
is feasible.

The crucial assumption that $A$ is essentially a low rank matrix
holds in many applications. For example, this is a model
hypothesis in the Latent Semantic Indexing (see \cite{BL, PRTV,
BDO, DDFLH, BDJ, AFKMS}). There $A$ is the ``document-term
matrix'', which is formed of the frequencies of occurrence of
various terms in the documents of a large collection. The
hypothesis that the documents are related to a small number of
(unknown) topics translates into the assumption that $A$ can be
approximated by an (unknown) low rank matrix. Finding such an
approximation would determine the ``best'' topics the collection
is really about. Other examples where this problem arises include
clustering of graphs \cite{DFKVV}, DNA microarray data, facial
recognition, web search (see \cite{DKM II}), lossy data
compression and cryptography (see \cite{BDJ}).

The best fixed rank approximation to $A$ is obviously given by the
partial sums of the Singular Value Decomposition (SVD)
$$
A = \sum_j \s_j(A) \; u_j \otimes v_j
$$
where $\s_j(A)$ is the nonincreasing and nonnegative sequence of
the singular values of $A$, and $u_j$ and $v_j$ are left and right
singular vectors of $A$ respectively. The best rank $k$ approximation
to $A$ in both the spectral and Frobenius norms is thus $A P_k$,
where $P_k$ is the orthogonal
projection onto the top $k$ left singular vectors of $A$.
In particular, for the spectral norm we have
\begin{equation}                    \label{best approx}
\min_{B: \, \rank B \le k} \|A - B\|_2 = \|A - A P_k\|_2 = \s_{k+1}(A).
\end{equation}

However, computing $P_k$, which gives the first elements of the
SVD of a $m \times n$ matrix $A$, is often impossible in practice
because (1) it would take many passes through $A$, which is
prohibitively slow for a matrix stored in an external memory; (2)
this would take superlinear time in $m+n$. Instead, it was
proposed in \cite{FKV, DK-conf, DKM II, DKM III} to use the
Monte-Carlo methodology: namely, approximate the $k$-th partial
sum of the SVD of $A$ by the $k$-th partial sum of the SVD of a
random submatrix of $A$. In this paper, we show that this can be
done:
\begin{enumerate}
\item with almost linear sample complexity $O(r \log r)$, that is
  by sampling only $O(r \log r)$ random rows of $A$, if $A$
  is approximable by a rank $r$ matrix;
\item in one pass through $A$ if the matrix is stored row-by-row,
  and in two passes if its entries are stored in arbitrary order;
\item using RAM space and time $O(n+m)$
  (and polynomial in $r$ and $k$).
\end{enumerate}

\begin{theorem}                     \label{t: low rank}
  Let $A$ be an $m \times n$ matrix with numerical rank
  $r = \norm{A}_F^2 / \norm{A}_2^2$.
  Let $\e, \d \in (0,1)$, and let $d \le m$ be an integer such that
  \begin{equation}              \label{i: complexity}
    d \ge  C \Bigl( \frac{r}{\e^4\d} \Bigr)
    \log \Bigl( \frac{r}{\e^4\d} \Bigr).
  \end{equation}
  Consider a $d \times n$ matrix $\tilde{A}$, which consists
  of $d$ normalized rows of $A$
  picked independently with replacement,
  with probabilities proportional to the squares of their Euclidean lengths.
  Then with probability at least $1 - 2\exp(-c/\d)$ the following holds.
  For a positive integer $k$, let $P_k$ be the orthogonal
  projection onto the top $k$ left singular vectors of $\tilde{A}$.
  Then
  \begin{equation}              \label{i: approx}
    \|A - AP_k\|_2 \le \s_{k+1}(A) + \e \norm{A}_2.
  \end{equation}
\end{theorem}
Here and in the sequel, $C, c, C_1, \ldots$ denote positive
absolute constants.

\medskip

Comparing \eqref{i: approx} with the best approximation
\eqref{best approx} given by the SVD, we see an additional error
$\e \norm{A}_2$ which can be made small by increasing the size $d$
of the sample.

\begin{remark}[Optimality]
  The almost linear sample complexity $d = O(r \log r)$
  achieved in Theorem~\ref{t: low rank} is optimal,
  see Proposition~\ref{optimality} below.
  The best known previous result, due to Drineas, Kannan and Mahoney,
  had with the quadratic sample complexity $d = O(r^2)$
  (\cite{DKM II} Theorem 4, see also \cite{DKM III} Theorem 3).
  The approximation scheme in Theorem \ref{t: low rank}
  was developed in \cite{FKV, DK-conf, DKM II, DKM III}.

\end{remark}

\begin{remark}[Numerical rank]
  The numerical rank $r = r(A) = \norm{A}_F^2 / \norm{A}_2^2$ in
  Theorem~\ref{t: low rank} is a relaxation of the exact notion of rank.
  Indeed, one always has $r(A) \le \rank(A)$. But as opposed to
  the exact rank, the numerical rank is stable under small perturbations
  of the matrix $A$. In particular, the numerical rank of $A$ tends to be low when
  $A$ is close to a low rank matrix, or when $A$ is sufficiently sparse.
  So results like Theorem~\ref{t: low rank}, which depend on the
  numerical rather than exact rank, should be useful in many applications,
  such as the Principal Component Analysis.
\end{remark}

\begin{remark}[Law of large numbers for operator-valued random variables]
  The new feature in our proof of Theorem \ref{t: low rank} is a use
  of the first author's argument about random vectors in the
  isotropic position \cite{R}. It yields a law of large numbers for
  operator-valued random variables. We apply it for independent
  copies of a rank one random operator, which is given by a random
  row of the matrix $A^T A$.
\end{remark}

\subsection{The cut-norm: decay, MAX-CSP problems}

Alon, Fernandez de la Vega, Kannan and Karpinski \cite{AFKK-conf,
AFKK-full} reduced the problem of additive approximation of the
MAX-CSP problems (which are NP-hard) to computing the cut-norm of
random submatrices. The cut norm of an $n \times n$ matrix $A$ is
the maximum sum of the entries of its submatrix,
$$
\|A\|_C = \max_{I, J} \Bigl| \sum_{i \in I, j \in J} A_{ij} \Bigr|,
$$
and it is equivalent to the $\ell_\infty \to \ell_1$ operator norm.

The problem is to understand how the cut norm of $A$ decreases when
we pass to its random submatrix. Let $Q$ be a a random subset of
$\{1,\ldots,n\}$ of expected cardinality $q$. This means that each
element of $\{1,\ldots,n\}$ is included into $Q$ independently with
probability $q/n$. We form a $Q \times Q$ random submatrix $A|_{Q
\times Q} = (A_{ij})_{i,j \in Q}$.

Intuitively, $A|_{Q \times Q}$ is $(q/n)^2$ times smaller
than $A$ if $A$ is diagonal-free, but only $(q/n)$ times smaller
than $A$ if $A$ is a diagonal matrix. We prove a general estimate of the
cut-norm of random submatrices, which combines both of these types of decay:

\begin{theorem}                     \label{t: CUT}
  Let $A$ be an $n \times n$ matrix.
  Let $Q$ be a random subset of $\{1,\ldots, n\}$
  of expected cardinality $q$.
  Then
  $$
  \E\|A|_{Q \times Q}\|_C
  \le O \Big(
      \Bigl( \frac{q}{n} \Bigr)^2 \|A - D(A)\|_C
    + \Bigl( \frac{q}{n} \Bigr) \|D(A)\|_C
    + \Bigl( \frac{q}{n} \Bigr)^{3/2} (\|A\|_\Col + \|A^T\|_\Col)
    \Big),
  $$
  where $\|A\|_\Col$ is the sum of the Euclidean lengths of the columns
  of $A$, and $D(A)$ is the diagonal part of $A$.
\end{theorem}

\begin{remark}[Optimality]
  The estimate in this theorem is optimal, see Section~\ref{s: cut optimality}.
\end{remark}

We now state a partial case of Theorem~\ref{t: CUT} in the form
useful for MAX-CSP problems. Note that $\norm{A}_{\Col} \le
\sqrt{n} \norm{A}_F$. Then we have:

\begin{corollary}                           \label{CSP}
  Under the hypotheses of Theorem~\ref{t: CUT}, let
  $q = \Omega(\e^{-2})$. Assume that
  $\|A\|_C = O(\e n^2)$, and $\|A\|_F = O(n)$, $\|A\|_\infty = O(\e^{-1})$,
  where $\|A\|_\infty$ denotes the maximum of the absolute
  values of the entries of $A$. Then
  $$
  \E\|A|_{Q \times Q}\|_C = O(\e q^2).
  $$
\end{corollary}

In solving MAX-2-CSP problems, one approximates the edge-weight
matrix $W$ of the graph on $n$ vertices by a cut approximation $W'$,
and checks that the the error matrix $A = W-W'$ satisfies the
assumptions of Corollary~\ref{CSP}, see \cite{F, AFKK-conf, AFKK-full}.
Hence the Corollary says that for a random induced graph on
$q = \Omega(\e^{-2})$ vertices, the same cut-approximation
(induced on the $q$ vertices of the random subgraph) works.
Namely, the error in cut-norm is at most $\e q^2$.

A weaker form of Corollary~\ref{CSP} was proved by Alon, Fernandez
de la Vega, Kannan and Karpinski \cite{AFKK-full}, Theorem 6.
Their result has a bigger sample complexity $q = \Omega(\e^{-4} \log(1/\e))$
and an extra assumption $n = e^{\Omega(\e^{-2})}$, but it works
for multidimensional arrays rather than for matrices (=$2$-dimensional
arrays).

Using Corollary~\ref{CSP} instead of \cite{AFKK-full} Theorem 6
slightly improves the best known sample complexity
for the approximation algorithm for MAX-2-CSP problems due to \cite{AFKK-full}.
The solution to a MAX-2SCP problem on $n$ variables can be approximated
within the additive error $\e n^2$ by the solution of the problem
induced by a randomly chosen $q$ variables.
The best known sample complexity, due to \cite{AFKK-full}, is
$q = O(\e^{-4} \log(1/\e))$. Using Corollary~\ref{CSP}
in the argument of \cite{AFKK-full} Theorem 6 improves the
sample complexity to $q = O(\e^{-4})$.

Our proof of Theorem~\ref{t: CUT} uses the technique of probability in
Banach spaces, and includes decoupling, symmetrization, and
application of a version of Slepian's lemma for Bernoulli random
variables due to Talagrand.

\subsection{The spectral norm: decay}

Perhaps the most important matrix norm is the spectral norm.
Nevertheless, its decay under
passing to submatrices has not been sufficiently understood.

Let $A$ be an $n \times n$ matrix, and $Q$ be a random subset of $\{1,\ldots,n\}$
of expected cardinality $q$ (as above).
We consider a random row-submatrix $A|_Q = (A_{ij})_{i \in Q, j \le n}$,
which consists of the rows of $A$ in $Q$.

When one orthogonally projects a vector $x \in \R^n$ onto $\R^Q$,
its Euclidean length reduces in average by the factor of
$\sqrt{\frac{q}{n}}$. So, one should expect a similar type of decay for the
spectral norm -- something like $\E\|A|_Q\|_2 \le \sqrt{\frac{q}{n}}
\|A\|_2$.

However, similarly to the previous section, the diagonal matrices
exhibit a different type of decay. For example, there is no decay at
all for the identity matrix. One can check that
the correct order of decay for diagonal matrices is
$$
\|A\|_{(k)} = \text{the average of $k$ biggest Euclidean lengths
of the columns of $A$,}
$$
where $k = n/q$. General matrices again
combine both types of decay:

\begin{theorem}  \label{t: norm decay}
  Let $A$ be an $n \times n$ matrix.
  Let $Q$ be a random subset of $\{1,\ldots, n\}$ of expected cardinality $q$.
  Then
  $$
  \E\|A|_Q\|_2
  \le O \Big( \sqrt{\frac{q}{n}} \|A\|_2
    + \sqrt{\log q} \, \|A\|_{(n/q)} \Big).
  $$
\end{theorem}

\begin{remark}[Optimality]
The estimate in this theorem is optimal. The example considered in
the proof of Proposition \ref{optimality} below shows that the
coefficient $\sqrt{\log q}$ is necessary.
\end{remark}

Generalizing an earlier result of Lunin \cite{L}, Kashin and
Tzafriri \cite{KT} (see \cite{V}) essentially proved the {\em
existence} of a subset $Q$ of cardinality $q$ and such that
$$
\|A|_Q\|_2 \le O \Big( \sqrt{\frac{q}{n}} \|A\|_2
    + \frac{\|A\|_F}{\sqrt{n}} \Big).
$$
Note that $\frac{\|A\|_F}{\sqrt{n}} = \big( \frac{1}{n}
\sum_{i=1}^n |A_i|^2 \big)^{1/2}$ is the average of the lengths of
{\em all} columns of $A$. As the example of diagonal operators
shows, for random subsets $Q$ this term has to be replaced by the
average of the few biggest columns. Talagrand \cite{T 95} proved
deep results on the more general operator norms $\ell_2 \to X$, where
$X$ is a $2$-smooth Banach space. However, the decay on
$\frac{q}{n}$ in his results is logarithmic rather than polynomial.

\subsection{Stable matrix results}

Many problems on random submatrices, of both theoretical and
practical importance, have functional-analytic rather than
linear-algeb\-raic nature. These problems, like those this paper
considers, are about estimating operator norms. We thus see a
matrix $A$ as a linear operator $A$ between finite dimensional
normed spaces -- say, between $\ell_2^n$ and $\ell_2^n$ for the spectral
norm, and between $\ell_\infty^n$ and $\ell_1^n$ for the cut norm.

From this perspective, the dimension $n$ of the ambient normed space
should play a minor role, while the real control of the picture
should be held by (hopefully few) quantities tied to the operator
rather than the space. As a trivial example, if $A$ is not of full
rank then the dimension $n$ is useless compared to the rank of $A$.
Further, we are looking for stable results, those not ruined by
small perturbations of the linear operators. This is a natural
demand in applications, and this differs our analytic perspective
from the linear algebraic one. It would thus be natural to look for
{\em stable} quantities tied to linear operators, which govern the
picture. For example, operator norms are stable quantities, while
the rank is not.

This paper advances such approach to matrices. The low
rank approximations in Theorem \ref{t: low rank} are only
controlled by the numerical rank $r(A) = \|A\|_F^2/\|A\|_2^2$
of the matrix, which is a stable relaxation of the rank.
The norms of random matrices in Theorems \ref{t: CUT} and
\ref{t: norm decay} are essentially controlled by the norms
of the original matrix (and naturally by the sampling factor,
the ratio of the size of the submatrix to the
size of the original matrix). The dimension $n$ of the matrix
does not play a separate role in these results (although
the matrix norms may grow with the dimension).

\medskip

{\bf Acknowledgement.}
This project started when the authors participated in
the PIMS Thematic Programme on Asymptotic Geometric Analysis
at the University of British Columbia in Summer 2002.
The first author was a PIMS postdoctoral fellow at that time.
We are grateful to PIMS for its hospitality.
The final part of this research was done when the first
authour visited University of California, Davis.
We are grateful for R.~Kannan for his comments
on the initial version of this paper, and to M.~Karpinski for
explaining what was the correct consequence of Corollary~\ref{CSP} for
MAX-2-CSP problems.
Finally, we thank the referees for their valuable comments
and suggestions.

\section{Notation}

For $p \le \infty$, the finite dimensional $\ell_p$ spaces
are denoted by $\ell_p^n$. Thus $\ell_p^n$ is the Banach space
$(\R^n, \|\cdot\|_p)$,
where $\|x\|_p = ( \sum_{i=1}^n |x_i|^p )^{1/p}$ for $p \le \infty$,
and $\|x\|_\infty = \max_i |x_i|$.
The closed unit ball of $\ell_p$ is denoted by
$B_p^n := \{x \,|\, \|x\|_p \le 1\}$.

The canonical basis of $\R^n$ is denoted by $(e_1,\ldots,e_n)$. Let
$x, y \in \R^n$. The canonical inner product is denoted by $\< x, y
\> := x^T y$. The tensor product is defined as $x \otimes y := y \,
x^T$; thus $(x \otimes y) z = \< x, z\> \, y$ for all $z \in \R^n$.

Let $A = (A_{ij})_{ij}$ be an $m \times n$ real matrix.
The spectral norm of $A$ is the operator norm $\ell_2 \to \ell_2$,
defined as
$$
\|A\|_2 := \sup_{x \in \R^n} \frac{\|Ax\|_2}{\|x\|_2} = \s_1(A),
$$
where $\s_1(A)$ is the largest singular value of $A$.
The Frobenius norm $\|A\|_F$ of $A$ is defined as
$$
\|A\|_F^2 := \sum_{i,j} A_{ij}^2 = \sum_j \s_j(A)^2,
$$
where $\s_j(A)$ are the singular values of $A$.

Finally, $C, C_1, c, c_1, \ldots$ denote positive absolute constants.
The $a = O(b)$ notation means that $a \le C b$ for some absolute constant $C$.

\section{Low rank approximations}            \label{s: lowrank}

In this section, we prove Theorem \ref{t: low rank},
discuss the algorithm for finding low rank approximations,
and show that the sample complexity in Theorem \ref{t: low rank} is optimal.
Our argument will be based on the law of large numbers
for operator-valued random variables.

\subsection{Law of large numbers for operator-valued random variables}
           \label{s: law of large numbers}

Theorem~\ref{t: low rank} is about random independent sampling the
rows of the matrix $A$. Such sampling can be viewed as an empirical
process taking values in the set of rows. If we sample enough rows,
then the matrix constructed from them would nicely approximate the
original matrix $A$ in the spectral norm. For the scalar random
variables, this effect is the classical Law of Large Numbers. For
example, let $X$ be a bounded random variable and let $X_1 \ldots
X_d$ be independent copies of $X$. Then
\begin{equation}                            \label{scalar LLN}
  \E \Big| \frac{1}{d} \sum_{j=1}^d X_j - \E X \Big|
  = O \Big( \frac{1}{\sqrt{d}} \Big).
\end{equation}
Furthermore, the large deviation theory allows one to estimate the
probability that the empirical mean $\frac{1}{d} \sum_{j=1}^d X_j$
stays close to the true mean $\E X$.

Operator-valued versions of this inequality are harder to prove.
The absolute value must be replaced by the operator norm.
So, instead of proving a large deviation estimate
for a single random variable, we have to estimate the
supremum of a random process. This requires deeper probabilistic techniques.
The following Theorem generalizes the main result of \cite{R}.

\begin{theorem}  \label{t: operator sampling}
  Let $y$ be a random vector in $\R^n$, which is uniformly bounded almost everywhere:
  $\|y\|_2 \le M$.
  Assume for normalization that $\norm{\E y \otimes y}_2 \le 1$.
  Let $y_1 \ldots y_d$ be independent copies of $y$. Let
  $$
  a := C \sqrt{\frac{\log d}{d}} \cdot M.
  $$
  Then

  (i) If $a < 1$ then
  \[
    \E \Big\| \frac{1}{d} \sum_{i=1}^d y_i \otimes y_i - \E \, y \otimes y \Big\|_2
    \le a.
  \]

  (ii) For every $t \in (0,1)$,
  $$
  \P \Big\{
    \Big\| \frac{1}{d} \sum_{i=1}^d y_i \otimes y_i - \E y \otimes y \Big\|_2
    > t \Big\}
  \le 2 \exp(-ct^2/a^2).
  $$
\end{theorem}

\begin{remark}
  Part (i) is a law of large numbers,
  and part (ii) is a large deviation estimate for operator-valued
  random variables. Comparing this result with its scalar prototype
  \eqref{scalar LLN},
  we see an additional logarithmic factor. This factor is essential,
  as we show in Remark~\ref{r: identity} below.
\end{remark}

\begin{remark}
  The boundedness assumption $\|y\|_2 \le M$ can be too strong for some
  applications. The proof of Theorem~\ref{t: operator sampling} shows that,
  in part (i), the boundedness almost everywhere can be relaxed to the moment
  assumption $\E \|y\|_2^q \le M^q$, where $q = \log d$.
  Part (ii) also holds under an assumption that the moments of $\|y\|_2$
  have a nice decay. However, we do not need these improvements here.
\end{remark}

\begin{remark} \label{r: identity}
The estimate in Theorem \ref{t: operator sampling} is in general
optimal.
Indeed, consider the random vector $y$ taking
values $\sqrt{n} e_1, \ldots, \sqrt{n} e_n$ each with probability
$1/n$, where $(e_i)$ is the canonical basis of $\R^n$.
Then $\E y \otimes y = I$. Then
\[
  \E \Big\| \frac{1}{d} \sum_{j=1}^d y_j \otimes y_j - I \Big\|_2
  = \E \max_{i=1 \ldots n}
  \left | \frac{n}{d} \,|\{ j \mid y_j=\sqrt{n} e_i\}| -1 \right |.
\]
If we want this quantity to be $O(1)$, then it is not hard to
check that $d$ should be of order at least $n \log n$. Therefore,
the coefficient $\sqrt{\log(d)/d}$ in Theorem \ref{t:
operator sampling} is optimal.
\end{remark}

\subsection{Proof of Theorem \ref{t: operator sampling}.}       \label{proof ops}
The proof consists of two steps. First we use the standard symmetrization
technique for random variables in Banach spaces, see e.g. \cite{LT} Section 6.
Then we adapt the technique of \cite{R} to obtain a bound on a symmetric
random process. To obtain the probability estimate in part (ii),
we shall estimate the high moments rather than the first moment
in part (i).

Let $\e_1 \ldots \e_d$ denote independent Bernoulli variables
taking values $1,-1$ with probability $1/2$.
Let $y_1 \ldots y_d, \bar y_1 \ldots \bar y_d$
be independent copies of $y$.
We shall denote by $\E_y$, $\E_{\bar{y}}$ and $\E_{\e}$
the expectations according to $(y_i)$, $(\bar{y}_i)$ and $(\e_i)$ respectively.

Let $p \ge 1$. We shall estimate
\begin{equation}                        \label{Ep}
E_p := \Big( \E \, \Big\| \frac{1}{d} \sumi \tens{y_i} - \E \tens{y} \Big\|_2^p
       \Big)^{1/p}.
\end{equation}
Note that
$
\E_y \, y \otimes y = \E_{\bar y} \, \bar y \otimes \bar y
= \E_{\bar y} \, \Big( \frac{1}{d} \sum_{i=1}^d \bar{y_i} \otimes \bar{y_i} \Big).
$
We put this into \eqref{Ep}. Since $x \mapsto \|x\|_2^p$ is a convex function
on $\R^n$, Jensen's inequality implies that
$$
E_p \le \Big( \E_y \E_{\bar y} \,
  \Big\| \frac{1}{d} \sumi \tens{y_i} - \frac{1}{d} \sumi \tens{\bar y_i}
  \Big\|_2^p \Big)^{1/p}.
$$
Since $\tens{y_i}-\tens{\bar y_i}$ is a symmetric random variable,
it is distributed identically with $\e_i(\tens{y_i}-\tens{\bar y_i})$.
Thus
$$
E_p \le \Big( \E_y \E_{\bar y} \E_{\e} \,
  \Big\| \frac{1}{d} \sumi \e_i(\tens{y_i}-\tens{\bar y_i})
  \Big\|_2^p \Big)^{1/p}.
$$
Denote $Y = \frac{1}{d} \sumi \e_i \tens{y_i}$
and $\bar Y = \frac{1}{d} \sumi \e_i \tens{\bar y_i}$.
Then
$\|Y - \bar Y\|_2^p
\le (\|Y\|_2 + \|\bar Y\|_2)^p
\le 2^p (\|Y\|_2^p + \|\bar Y\|_2^p)$,
and $\E \|Y\|_2^p = \E \|\bar Y\|_2^p$.
Thus we obtain
$$
E_p \le 2 \Big( \E_y \E_\e \; \Big\| \frac{1}{d} \sumi \e_i \tens{y_i}
  \Big\|_2^p \Big)^{1/p}.
$$

We shall estimate the last expectation using a lemma from \cite{R}.

\begin{lemma}  \label{l: Khinchine}
Let $y_1 \ldots y_d$ be vectors in $R^k$ and let $\e_1 \ldots \e_d$
be independent Berno\-ulli variables taking values $1,-1$ with
probability $1/2$. Then
\[
\Big( \E \Big\|\sum_{i=1}^d \e_i \tens{y_i} \Big\|_2^p \Big)^{1/p}
\le C_0 (p+\log k)^{1/2}
  \cdot \max_{i=1 \ldots d} \norm{y_i}_2
  \cdot \Big\| \sum_{i=1}^d \tens{y_i} \Big\|_2^{1/2}.
\]
\end{lemma}
\begin{remark}   \label{r: number of vectors}
 We can consider the vectors $y_1 \ldots y_d$
 as vectors in their linear span, so we can always choose
 the dimension $k$ of the ambient space at most $d$.
\end{remark}

Combining Lemma~\ref{l: Khinchine} with Remark \ref{r: number of
vectors} and using H\"older's inequality, we obtain
\begin{equation} \label{enorm}
E_p \le 2C_0 \frac{(p + \log d)^{1/2}}{d} \cdot M \cdot
  \Big( \E \Big\| \sumi \tens{y_i} \Big\|_2^{p} \Big)^{1/2p}.
\end{equation}
By Minkowski's inequality we have
$$
\Big( \E \Big\| \sumi \tens{y_i} \Big\|_2^{p} \Big)^{1/p}
\le d \Big[ \Big( \E \Big\| \frac{1}{d} \sumi \tens{y_i} - \E \, \tens{y}
  \Big\|_2^p \Big)^{1/p}
  + \| \E \, \tens{y} \|_2
  \Big]
\le d(E_p + 1).
$$
So we obtain
$$
E_p \le \frac{a p^{1/2}}{2} (E_p + 1), \qquad
\text{where} \quad
a = 4C_0 \Big( \frac{\log d}{d} \Big)^{1/2} M.
$$
It follows that
\begin{equation}                    \label{EpEp}
  \min(E_p,1) \le a p^{1/2}.
\end{equation}

To prove part (i) of the theorem, note that $a \le 1$ by the assumption.
It thus follows that $E_1 \le a$.
This proves part (i).

To prove part (ii), we can $E_p = (\E \, Z^p)^{1/p}$,
where
$$
Z = \Big\| \frac{1}{d} \sumi \tens{y_i} - \E \tens{y} \Big\|_2.
$$
So \eqref{EpEp} implies that
\begin{equation}                    \label{Ep 1}
  \big (\E  \min(Z,1)^p \big )^{1/p}
  \le \min (E_p, 1)
  \le ap^{1/2}.
\end{equation}
This moment bound can be expressed as a tail probability estimate
using the following standard lemma, see e.g. \cite{LT} Lemmas 3.7
and 4.10.

\begin{lemma}
  Let $Z$ be a nonnegative random variable.
  Assume that there exists a constant $K > 0$ such that
  $(\E \, Z^p)^{1/p} \le K p^{1/2}$ for all $p \ge 1$.
  Then
  $$
  \P \{Z > t\} \le 2 \exp(-c_1 t^2/K^2)
  \qquad \text{for all $t > 0$}.
  $$
\end{lemma}

It thus follows this and from \eqref{Ep 1} that
$$
\P \{ \min(Z,1) > t \} \le 2 \exp (-c_1 t^2/a^2)
\qquad \text{for all $t > 0$}.
$$
This completes the proof of the theorem.
\endproof

\subsection{Proof of Theorem \ref{t: low rank}}

By the homogeneity, we can assume $\|A\|_2 = 1$.

The following lemma of Drineas and Kannan \cite{DK-conf} (see also
\cite{DKM II}) reduces Theorem \ref{t: low rank} to a comparison
of $A$ and a sample $\tilde{A}$ in the spectral norm.

\begin{lemma}[Drineas, Kannan]                  \label{DK}
  $$
  \|A - A P_k\|_2^2
  \le \s_{k+1}(A)^2 + 2 \|A^T A - \tilde{A}^T \tilde{A}\|_2.
  $$
\end{lemma}

\begin{proof}
We have
\begin{align*}
\|A - A P_k\|_2^2
  &=  \sup_{x \in \ker P_k, \; \|x\|_2 = 1} \|Ax\|_2^2
   =  \sup_{x \in \ker P_k, \; \|x\|_2 = 1} \< A^T Ax, x \>  \\
  &\le  \sup_{x \in \ker P_k, \; \|x\|_2 = 1} \< (A^TA - \tilde{A}^T \tilde{A})x, x \>
     + \sup_{x \in \ker P_k, \; \|x\|_2 = 1} \< \tilde{A}^T \tilde{A} x, x \> \\
  &=  \|A^T A - \tilde{A}^T \tilde{A}\|_2 + \s_{k+1} (\tilde{A}^T \tilde{A}).
\end{align*}

By a result of perturbation theory, $| \s_{k+1} (A^T A) - \s_{k+1}
(\tilde{A}^T \tilde{A}) | \le  \|A^T A - \tilde{A}^T \tilde{A}\|_2$.
This proves Lemma~\ref{DK}.
\end{proof}

Let $x_1 \ldots x_m$ denote the rows of the matrix $A$.
Then
\[
   A^T A = \sum_{j=1}^m x_j \otimes x_j.
\]
We shall view the matrix $A^TA$ as the {\em true mean} of a bounded
operator valued random variable, whereas $\tilde{A}^T \tilde{A}$ will be
its {\em empirical mean}; then we shall apply the Law of Large Numbers
for operator-valued random variables -- Theorem~\ref{t: operator sampling}.
To this end, define a random vector $y \in \R^m$ as
\[
  \P \left ( y = \frac{\norm{A}_F}{\norm{x_j}_2} \, x_j \right )
  = \frac{\norm{x_j}_2^2}{\norm{A}_F^2}.
\]
Let $y_1 \ldots y_d$ be independent copies of $y$.
Let the matrix $\tilde{A}$ consist of rows
$\frac{1}{\sqrt{d}} y_1 \ldots \frac{1}{\sqrt{d}}  y_d$.
(The normalization of $\tilde{A}$ here is different than in the statement
of Theorem \ref{t: low rank}: in the proof, it is convenient to
multiply $\tilde{A}$ by the factor $\frac{1}{\sqrt{d}} \|A\|_F$.
However note that the singular vectors of $\tilde{A}$ and
thus $P_k$ do not change.)
Then
$$
A^T A = \E \tens{y}, \qquad
\tilde{A}^T \tilde{A} = \frac{1}{d} \sumi \tens{y_j}, \qquad
M := \norm{y}_2 = \norm{A}_F = \sqrt{r}.
$$
We can thus apply Theorem \ref{t: operator sampling}. Due to our
assumption on $d$, we have
$$
a := 4C_0 \Big( \frac{\log d}{d} \cdot r \Big)^{1/2} \le \frac{\e^2
\d^{1/2}}{2} < 1.
$$
Thus Theorem \ref{t: operator sampling} yields (with $t = \e^2/2$)
that, with probability at least $1 - 2 \exp(-c/\d)$, we have
\[
  \| \tilde{A}^T \tilde{A} - A^T A \|_2
  \le \frac{\e^2}{2}.
\]
Whenever this event holds, we can conclude by Lemma~\ref{DK} that
 \[
  \|A - A P_k\|_2
    \le \s_{k+1}(A) + \sqrt{2} \|A^T A - \tilde{A}^T \tilde{A}\|_2^{1/2}
    \le \s_{k+1}(A) + \e.
  \]
This proves Theorem \ref{t: low rank}.
\endproof

\subsection{Algorithmic aspects of Theorem \ref{t: low rank}.}
Finding a good low rank approximation to a matrix $A$
amounts, due to Theorem \ref{t: low rank}, to sampling a random
submatrix $\tilde{A}$ and computing its SVD (actually, only
left singular vectors are needed).
The algorithm works well if the numerical rank
$r = r(A) = \|A\|_F^2/\|A\|_2^2$ of the matrix $A$ is small.
This is the case, in particular, when $A$ is essentially a
low-rank matrix, because $r(A) \le \rank(A)$.

First, the algorithm samples $d = O(r \log r)$ random rows of $A$.
Namely, it takes $d$ independent samples of the random vector
$y$ whose law is
$$
\P \left ( y = \frac{A_j}{\|A_j\|_2} \right )
= \frac{\|A_j\|_2^2}{\norm{A}_F^2}
$$
where $A_j$ is the $j$-th row of $A$.
This sampling can be done in {\em one pass} through $A$
if the matrix is stored row-by-row, and in {\em two passes}
if its entries are stored in arbitrary order \cite[Section 5.1]{DFKVV}.

Then the algorithm computes the SVD of the $d \times n$
matrix $\tilde{A}$, which consists of the normalized sampled rows.
This can be done in time $O(dn) + $ the time needed to compute the
SVD of a $d \times d$ matrix. The latter can be done by one of the known methods.
This takes significantly less time than computing
SVD of the original $m \times n$ matrix $A$.
In particular, the running time of this algorithm
is linear in the dimensions of the matrix (and polynomial
in $d$).

\subsection{Optimality of the sample complexity}

The sample complexity $d = O(r \log r)$ in
Theorem~\ref{t: low rank} is best possible:

\begin{proposition}                         \label{optimality}
  There exist matrices $A$ with arbitrarily big numerical
  rank $r = \|A\|_F^2 / \|A\|_2^2$ and such that whenever
  $$
  d < \frac{1}{10} r \log r,
  $$
  the conclusion \eqref{i: approx} of Theorem~\ref{t: low rank}
  fails for $k=n$ and for all $\e \in (0,1)$.
\end{proposition}

\medskip

\proof Let $n, m \in \N$ be arbitrary numbers such that $n<m$. We
define the $m \times n$ matrix  by its entries as follows:
$$
A_{ij} = \sqrt{\frac{n}{m}} \; \delta_{\lceil \frac{n}{m}i \rceil,
          j},
$$
where $\d_{ij} = 1$ if $i = j$ and $\d_{ij} = 0$ otherwise.

Then each row of $A$ contains exactly one entry of value $\sqrt{\frac{n}{m}}$,
and each row is repeated $m/n$
times. The $j$-th column of $A$ contains exactly one block of values
$\sqrt{\frac{n}{m}}$ in positions
$i \in (\frac{m}{n}(j-1), \frac{m}{n}\;j ] =: I_j$.
In particular, the columns are orthonormal.
Also, $\|A\|_2 = 1$, $\|A\|_F = \sqrt{n}$, thus $r = n$.

Now we form a submatrix $\tilde{A}$ as described in Theorem 1.1 --
by picking $d$ rows of $A$ independently and with uniform distribution.
If $d < \frac{1}{10} n \log n$, then with high probability there exists
at least one block $I_j$ from which no rows $i$ are picked. Call this block
$I_{j_0}$. It follows that $j_0$-th column of $\tilde{A}$ is zero.
Consider the coordinate vector $e_{j_0} = (0,\ldots,0,1,0,\ldots,0)$
of $n$ positions, with $1$ at position $j_0$.
Then
$e_{j_0} \in \ker \tilde{A} \subseteq \ker P_k \subseteq \ker (AP_k)$.
Thus
$\|(A - AP_k)e_{j_0}\|_2 = \|Ae_{j_0}\|_2 = 1$. Hence
$$
\|A - AP_k\|_2 \ge 1,
\ \ \ \text{while} \ \ \
\s_{n+1}(A) = 0, \ \
\|A\|_2 = 1.
$$
Hence \eqref{i: approx} fails for $k=n$ and for all $\e \in (0,1)$.
\endproof

\section{The decay of the cut norm}             \label{s: cut}

In this section, we prove Theorem \ref{t: CUT} on the cut norm
of random submatrices and show that it is optimal.
Our argument will be based on the tools of probability in
Banach spaces: decoupling, symmetrization, and Slepian's Lemma
(more precisely, its version for the Rademacher random variables
due to M.Talagrand).

\subsection{Proof of Theorem \ref{t: CUT}}

It is known and easy to check that
$$
\frac{1}{4}\|A\|_{\infty \to 1} \le \|A\|_C \le \|A\|_{\infty \to 1},
$$
where $\|A\|_{\infty \to 1}$ denotes the operator norm of $A$
from $\ell_\infty^n$ into $\ell_\infty^n$:
$$
\|A\|_{\infty \to 1} := \sup_{x \in \R^n} \frac{\|Ax\|_1}{\|x\|_\infty}
= \sup_{x \in B_\infty^n} \|Ax\|_1
$$
(recall that $B_\infty^n$ denotes the unit ball of $\ell_\infty^n$).
Note also that both these norms are self-dual:
$$
\|A^T\|_C = \|A\|_C, \qquad \|A^T\|_{\infty \to 1} = \|A\|_{\infty \to 1}.
$$
So we can prove Theorem~\ref{t: CUT} for the norm $\|\cdot\|_{\infty \to 1}$
instead of the cut norm.

We shall use the following decoupling lemma due to Bourgain and Tzafriri \cite{BT}.

\begin{lemma}  \label{l: decoupling}
  Let $(\xi_i)$ be a finite sequence of bounded i.i.d. random variables,
  and $(\xi'_i)$ be its independent copy.
  Then for any sequence of vectors $(x_{ij})$ in a Banach space
  with $x_{ii} = 0$,
  $$
  \E \Big\| \sum_{i,j} \xi_i \xi_j x_{ij} \Big\|
  \le  20 \E \Big\| \sum_{i,j} \xi_i \xi'_j x_{ij} \Big\|.
  $$
\end{lemma}

Let $\d_1 \ldots \d_n$ be independent Bernoulli random variables,
which take value $1$ with probability $\d:=q/n$.
Let $P_\D$ denote the coordinate projection on the random
set of coordinates $\{j \mid \d_j=1\}$.

Denote by $D(A)$ the diagonal part of $A$. Then
$$
P_\D A P_\D = P_\D (A - D(A))P_\D + P_\D D(A) P_\D = \sum_{i \ne
j} \d_i \d_j A_{ij} e_i \otimes e_j
  + \sum_{i=1}^n \d_i A_{ii} e_i \otimes e_i.
$$
We can use Lemma \ref{l: decoupling} to estimate the first
summand, taking $x_{ij} = A_{ij} e_i \otimes e_j$ if $i \ne j$ and
$x_{ij} = 0$ if $i = j$.
To this end, let $(\d'_j)$ be an independent copy of $(\d_j)$,
and let $P_{\D'}$ denote the coordinate projection on the random
set of coordinates $\{j \mid \d'_j=1\}$.
Then by Lemma \ref{l: decoupling} and by the triangle inequality we obtain
$$
\E \|P_\D A P_\D\|_{\infty \to 1} \le 20 \E \|P_\D (A - D(A))
P_{\D'}\|_{\infty \to 1}
    + \d \sum_{i=1}^n |A_{ii}|.
$$
Clearly, $\sum_{i=1}^n |A_{ii}| = \|D(A)\|_{\infty \to 1}$. Thus to
complete the proof, we can assume that the diagonal of $A$ is zero,
and prove the inequality as stated in the theorem for $\E \| P_\D A
P_{\D'} \|_{\infty \to 1}$, i.e.
\begin{equation}                        \label{diag free}
  \E \| P_\D A P_{\D'} \|_{\infty \to 1}
  \le C \d^2 \|A\|_{\infty \to 1} + C \d^{3/2} (\|A\|_\Col + \|A^T\|_\Col).
\end{equation}
Note that
$$
\E \| A P_{\D'} \|_{\infty \to 1}
= \E \sup_{x \in  B_\infty^n} \sum_{i=1}^n |\< A P_{\D'} x, e_i \> |,
$$
hence
\begin{align}                           \label{PAP}
\E \|P_\D A P_{\D'}\|_{\infty \to 1}
  &=  \E \sup_{x \in  B_\infty^n}
         \sum_{i=1}^n \d_i |\< A P_{\D'} x, e_i \> |  \notag\\
  &=  \E \sup_{x \in  B_\infty^n}
          \sum_{i=1}^n (\d_i - \d) |\< A P_{\D'} x, e_i \> |
      + \d \cdot \E \|A P_{\D'}\|_{\infty \to 1}.
\end{align}
We proceed with a known symmetrization argument, which
we used in the beginning of Section~\ref{proof ops}.
Since $\d_i - \d$ are mean zero, we can replace $\d$ by $\d''_i$,
an independent copy of $\d_i$, which can only increase the quantity
in \eqref{PAP}. Then the first term in \eqref{PAP} does not exceed
\begin{equation}                    \label{EED}
\E \sup_{x \in  B_\infty^n}
          \sum_{i=1}^n (\d_i - \d_i'') |\< A P_{\D'} x, e_i \> |.
\end{equation}
The random variable  $\d_i - \d''_i$ is symmetric, hence it is
distributed identically with $\e_i(\d_i - \d''_i)$, where $\e_i$
are $-1,1$-valued symmetric random variables independent of
all other random variables. Therefore the expression in \eqref{EED}
bounded by
\begin{multline}                \label{EEE}
\E \sup_{x \in  B_\infty^n}
          \sum_{i=1}^n \e_i \d_i |\< A P_{\D'} x, e_i \> |
  + \E \sup_{x \in  B_\infty^n}
          \sum_{i=1}^n \e_i \d''_i |\< A P_{\D'} x, e_i \> | \\
\le 2 \E \sup_{x \in  B_\infty^n} \sum_{i=1}^n \e_i \d_i |\< A P_{\D'} x, e_i \> |.
\end{multline}
To estimate this, we use Slepian's inequality for Rademacher random
variables proved by Talagrand. This estimate allows us to remove the
absolute values in \eqref{EEE}. Precisely, a partial case of
Slepian's inequality due to Talagrand (see \cite{LT}, equation
(4.20)) states that, for arbitrary $y_1, \ldots, y_n \in \R^n$, one
has
$$
\E \sup_{x \in B_\infty^n} \sum_{i=1}^n \e_i |\< x, y_i\> |
\le  \E \sup_{x \in B_\infty^n} \sum_{i=1}^n \e_i \< x, y_i\>
=  \E \Big\| \sum_{i=1}^n \e_i y_i \Big\|_1.
$$
Therefore
\begin{align*}
\E \sup_{x \in  B_\infty^n} \sum_{i=1}^n \e_i \d_i |\< A P_{\D'} x, e_i \> |
  &=   \E \sup_{x \in  B_\infty^n} \sum_{i=1}^n \e_i \Big| \< x, P_{\D'} A^T \d_i e_i \> \Big| \\
  &\le \E \Big\| P_{\D'} A^T \Big( \sum_{i=1}^n \e_i \d_i e_i \Big) \Big\|_1 \\
  &=   \E \sum_{j=1}^n \d'_j \Big| \< A^T (\sum_{i=1}^n \e_i \d_i e_i), e_j \> \Big|\\
  &=   \d \cdot \E \sum_{j=1}^n \Big| \sum_{i=1}^n \e_i \d_i A_{ij} \Big| \\
  &\le \d \cdot \sum_{j=1}^n \Big( \E \Big| \sum_{i=1}^n \e_i \d_i A_{ij} \Big|^2 \Big)^{1/2} \\
  &=   \d \cdot \sum_{j=1}^n \Big( \E \sum_{i=1}^n |\d_i A_{ij}|^2 \Big)^{1/2}
    \quad \text{(averaging over $(\e_i)$)}\\
  &=   \d^{3/2} \cdot \sum_{j=1}^n \Big( \sum_{i=1}^n |A_{ij}|^2 \Big)^{1/2}
  = \d^{3/2} \|A\|_\Col.
\end{align*}
We have proved that the first term in \eqref{PAP} does not exceed
$\d^{3/2} \|A\|_\Col$.
To estimate the second term, note that
$$
\E \|A P_{\D'}\|_{\infty \to 1}
= \E \|P_\D A^T\|_{\infty \to 1}
= \E \sup_{x \in  B_\infty^n} \sum_{i=1}^n \d_i | \< A^T x, e_i \> |.
$$
So we can essentially repeat the argument above to bound this expression by
$$
\le \d^{1/2} \|A^T\|_\Col + \d \|A^T\|_{\infty \to 1}
=   \d^{1/2} \|A^T\|_\Col + \d \|A\|_{\infty \to 1}.
$$
Putting this together, we can estimate \eqref{PAP} as
\begin{multline*}
\E \|P_\D A P_{\D'}\|_{\infty \to 1}
  \le \d^{3/2} \|A\|_\Col + \d (\d^{1/2} \|A^T\|_\Col + \d \|A\|_{\infty \to 1}) \\
  \le \d^{3/2} \|A\|_\Col + \d^{3/2} \|A^T\|_\Col + \d^2 \|A\|_{\infty \to 1},
\end{multline*}
as desired.
This completes the proof of Theorem~\ref{t: CUT}.
\endproof

\subsection{Optimality}         \label{s: cut optimality}
All terms appearing in Theorem \ref{t: norm decay} are necessary.
Their optimality can be witnessed on different types of matrices. To
see that the first term is necessary, consider a matrix $A$, all
whose entries are equal $1$. For this matrix $\norm{A}_C=n^2$, and
for any $Q \subset \{1, \ldots n\}, \ \norm{A_{Q \times
Q}}_C=|Q|^2$.

The optimality of the second term can be seen in the case when $A$
is the identity matrix. In this case $\norm{A}_C=n$, while
$\norm{A_{Q \times Q}}=|Q|$.

To prove that the third term is also necessary, assume that
$A=(\e_{i,j})$ is a random $\pm 1$ matrix. Then $\norm{D(A)}_C=n$,
and $\norm{A}_{\Col} =\norm{A^T}_{\Col} =n^{3/2}$. It is easy to
show that $\E_{\e}\norm{A}_C \le C n^{3/2}$, so for $q<n$ the third
term in Theorem \ref{t: norm decay} is dominant. Indeed, by Azuma's
inequality, for any $x,y \in \{0,1\}^n$
\[
  \P_{\e} \left ( \left | \frac{1}{n}
             \sum_{i,j=1}^n \e_{ij}x_i y_j \right | >t
          \right )
  \le C e^{-t^2/2}.
\]
Hence,
\[
  \P_{\e} \big ( \norm{A}_C > s n^{3/2} \big )
  \le 4^n \cdot C e^{-s^2 n},
\]
which implies the desired bound for the expectation.

Now fix a $\pm 1$ matrix $A$ such that $\norm{A}_C \le C n^{3/2}$.
Let $Q$ be any subset of $\{1, \ldots, n\}$. Recall that the norms
$\norm{A}_C$ and $\norm{A}_{\infty \to 1}$ are equivalent. We claim
that
\[
  \norm{A|_{Q \times Q}}_{\infty \to 1}
  \ge \frac{1}{\sqrt{2}} |Q|^{3/2}.
\]
Indeed, let $\d_i, \ i \in Q$ be independent $\pm 1$ random
variables. Then by Khinchine's inequality
\[
  \sum_{j \in Q} \E_{\d} \left | \sum_{i \in Q} \e_{ij} \d_i
                         \right |
  \ge \frac{1}{\sqrt{2}} |Q|^{3/2}.
\]
Choose $x \in \{-1,1\}^Q$ such that $\sum_{j \in Q} \left | \sum_{i
\in Q} \e_{ij} x_i \right | \ge \frac{1}{\sqrt{2}} |Q|^{3/2}$. For
$j \in Q$ set
\[
  y_j= \text{sign} \left ( \sum_{i \in Q} \e_{ij} x_i \right ).
\]
Then
\[
  \norm{A|_{Q \times Q}}_{\infty \to 1}
  \ge \left | \sum_{i,j \in Q} \e_{ij} x_i y_j \right |
  \ge \frac{1}{\sqrt{2}} |Q|^{3/2}.
\]
Therefore,
\[
  \E_Q \norm{A|_{Q \times Q}}_C
  \ge \frac{1}{4 \sqrt{2}} \left ( \frac{q}{n} \right )^{3/2} \cdot
      \big ( \norm{A}_{\Col} + \norm{A^T}_{\Col} \big ).
\]
Therefore, the third term is also necessary.

\section{The decay of the spectral norm}

In this section, we prove Theorem \ref{t: norm decay} on the
spectral norm of random submatrices.

By homogeneity we can assume that $\norm{A}_2=1$.
Let $\d_1,\ldots,\d_n$ be $\{ 0, 1 \}$-valued independent random variables
with $\E \d_j = \d = \frac{q}{n}$. So our random set is
$Q = \{j \, | \, \d_j = 1\}$.

Let $x_1 \ldots x_n$ denote the columns of $A$. Then
\[
  A  =  \sum_{j=1}^n e_j \otimes x_j, \qquad
  A|_Q  =  \sum_{j=1}^n \d_j e_j \otimes x_j.
\]
The spectral norm can be computed as
$$
\|A\|_2 = \|A^T A\|_2^{1/2} = \Big\| \sum_{j=1}^n x_j \otimes x_j \Big\|_2^{1/2},
$$
and similarly
$$
\|A|_Q\|_2 = \Big\| \sum_{j=1}^n \d_j x_j \otimes x_j \Big\|_2^{1/2}.
$$
To estimate the latter norm, we shall first apply the standard
symmetrization argument (see \cite{LT} Lemma 6.3), like we did in
the beginning of Section~\ref{proof ops} and in Section \ref{s:
cut}. Then we will apply Lemma \ref{l: Khinchine}. Set
\[
  E=\E \|A|_Q\|_2.
\]
The symmetrization argument yields
$$
E \le  \E \Big\| \sum_{j=1}^n
                  (\d_j - \d) x_j \otimes x_j \Big\|_2^{1/2}
        +  \sqrt{\d} \norm{A}_2^{1/2}
  \le  2 \E_\d \Big( \E_\e
          \Big\| \sum_{j=1}^n \e_j \d_j x_j \otimes x_j
          \Big\|_2 \Big)^{1/2}
        +  \sqrt{\d}.
$$

Now we apply Lemma \ref{l: Khinchine} with $p=1$ to bound
$\E_\e \Big\| \sum_{j=1}^n \e_j \d_j x_j \otimes x_j \Big\|_2$
for fixed $(\d_j)$. By Remark \ref{r: number of vectors}, we can assume
$k$ in this Lemma equal
$$
n(\d) :=  e+  \sum_{j \le n} \d_j.
$$
Then using Cauchy-Schwartz  inequality we obtain
\begin{align}                     \label{eless}
  E &\le \E_\d \Big( C \sqrt{\log n(\d)} \cdot
                      \max_{j=1 \ldots n} \d_j \norm{x_j}_2 \cdot
                      \Big\| \sum_{j=1}^n \d_j x_j \otimes  x_j \Big\|_2^{1/2}
              \Big)^{1/2}  + \sqrt{\d}  \notag \\
    &\le C \Big( \E_\d \Big( \sqrt{\log n(\d)} \cdot
            \max_{j=1 \ldots n} \d_j \norm{x_j}_2
           \Big)  \Big)^{1/2}
           \Big( \E_{\d} \Big\| \sum_{j=1}^n \d_j x_j \otimes x_j \Big\|_2^{1/2}
           \Big)^{1/2} + \sqrt{\d}.
\end{align}
To estimate the fist term in the product here, we use the following

\begin{lemma}   \label{l: order statistics}
Let $a_1 \ge a_2 \ge \ldots \ge a_n \ge 0$ and let $\d_1 \ldots
\d_n$ be independent Bernoulli random variables taking value $1$
with probability $\d>2/n$. Then
\[
      \frac{\d}{4e} \sqrt{\log \d n} \cdot \sum_{j=1}^{1/\d} a_j
  \le \E \left (\sqrt{\log n(\d)}
                \cdot \max_{j=1 \ldots n}  \d_j a_j \right )
  \le 4\d \sqrt{\log \d n} \cdot \sum_{j=1}^{1/\d} a_j.
\]
\end{lemma}
\begin{proof}
To prove the upper estimate note that
\[
  \max_{j=1 \ldots n} \d_j a_j
  \le
  \sum_{j=1}^{1/\d} \d_j a_j + a_{1/\d}.
\]
Hence,
\begin{equation}  \label{stat1}
   \E \left (\sqrt{\log n(\d)}
                \cdot \max_{j=1 \ldots n}  \d_j a_j \right )
   \le \E \Big(\sqrt{\log n(\d)}
                \cdot \sum_{j=1}^{1/\d} \d_j a_j \Big)
   + a_{1/\d} \cdot \E \sqrt{\log n(\d)}.
\end{equation}
Jensen's inequality yields
\begin{equation}                            \label{star}
  \E \sqrt{\log n(\d)}
  \le  \sqrt{\log (\E \sum_{i=1}^n \d_i+e)} \le 2 \sqrt{\log \d n}.
\end{equation}
By the linearity of expectation, the first term in the right hand side
of \eqref{stat1} equals
\[
  \sum_{j=1}^{1/\d} a_j
        \E \left (\d_j \sqrt{\log n(\d)} \right)
  \le \sum_{j=1}^{1/\d} a_j
        \E \left (\d_j \sqrt{\log (\sum_{i \neq j} \d_i +1+e)} \right),
\]
where we estimated $n(\d)$ replacing $\d_j$ by 1. Taking the
expectation first with respect to $\d_j$ and then with respect to
the other $\d_i$, and using Jensen's inequality,
we bound the last expression by
\begin{equation}                        \label{double-star}
  \d \sum_{j=1}^{1/\d} a_j \cdot \sqrt{\log (\d n+1+e)}
  \le 2 \d \sum_{j=1}^{1/\d} a_j \cdot \sqrt{\log \d n}.
\end{equation}
Finally, substituting \eqref{star} and \eqref{double-star} into \eqref{stat1},
we obtain
\[
   \E \left (\sqrt{\log n(\d)}
                \cdot \max_{j=1 \ldots n}  \d_j a_j \right )
   \le \Big( 2 \d \sum_{j=1}^{1/\d} a_j + 2a_{1/\d} \Big)
              \cdot \sqrt{\log \d n}
   \le 4 \d \sum_{j=1}^{1/\d} a_j \cdot \sqrt{\log \d n}.
\]

To prove the lower bound, we estimate the product in Lemma \ref{l:
order statistics} from below to make the terms independent. We
have
\begin{align} \label{stat2}
  \E \left (\sqrt{\log n(\d)}
                \cdot \max_{j=1 \ldots n}  \d_j a_j \right )
  &\ge \E \left (\sqrt{\log (\sum_{i=1/\d+1}^n \d_i+e)}
                \cdot \max_{j=1 \ldots 1/\d}  \d_j a_j \right )\notag \\
  &= \E \sqrt{\log (\sum_{i=1/\d+1}^n \d_i+e)}
                \cdot \E \max_{j=1 \ldots 1/\d}  \d_j a_j.
\end{align}
These terms  will be estimated separately. Since $\P \, (
\sum_{i=1/\d+1}^n \d_i \ge \d n/2) \ge 1/2$,
\[
  \E \sqrt{\log (\sum_{i=1/\d+1}^n \d_i+e)}
  \ge \frac{1}{2} \sqrt{\log \frac{\d n}{2}}.
\]

Let  $1 \le k \le 1/\d$. Denote by $A_k$ the event
$\{ \d_k=1, \d_j=0 \text{ for } 1 \le j \le 1/\d, \ j \neq k\}$. Then
\[
  \P (A_k) = \d \cdot (1-\d)^{1/\d-1} \ge \d/e.
\]
Since the events $A_1, \ldots, A_{1/\d}$ are disjoint,
\[
  \E \max_{j=1 \ldots 1/\d} \d_j a_j \ge
  \sum_{k=1}^{1/\d} a_k \P (A_k) \ge
  \frac{\d}{e} \sum_{j=1}^{1/\d} a_j.
\]
Substituting this estimate into \eqref{stat2} finishes the proof
of Lemma~\ref{l: order statistics}.
\end{proof}

Now we can complete the proof of Theorem \ref{t: norm decay}.
Combining Lemma \ref{l: order statistics} and \eqref{eless}, we
get
$$
E \le C \Big( 4\d \sqrt{\log \d n} \sum_{j=1}^{1/\d} \|x_j\|_2 \Big)^{1/2} E^{1/2} + \sqrt{\d}
= 2C \Big( \sqrt{\log \d n} \|A\|_{(1/\d)} \Big)^{1/2} E^{1/2} + \sqrt{\d}.
$$
It can be easily checked that $E \le a E^{1/2}+b$ implies $E \le
4a^2 +2b$. Hence, recalling that $\d=q/n$, we conclude that
\[
E \le 16 C^2 \sqrt{\log q} \cdot \norm{A}_{(n/q)} + 2\sqrt{q/n}.
\]
This completes the proof of Theorem \ref{t: norm decay}.
\endproof

\end{document}